\documentclass{amsart}
\usepackage{amssymb, amsmath, amscd}
\usepackage{latexsym}
\usepackage{graphicx}

\newcommand{\BC}{\Bbb C}

\newcommand{\supp}{\operatorname{supp}}

\newcommand{\repn}{representation\,}
\newcommand{\repns}{representations\,}

\newtheorem{prop}{Proposition}[section]
\newtheorem{theo}{Theorem}[section]
\newtheorem{lemm}{Lemma}[section]
\newtheorem{cor}{Corollary}[section]

\begin{document}

\title[restricted algebras]{the restricted algebras on inverse semigroups III, Fourier algebra}
\date{}
\author[M. Amini, A.R. Medghalchi]{Massoud Amini, Alireza Medghalchi}
\address{Department of Mathematics, Tarbiat Modarres University, P.O.Box 14115-175,
Tehran, Iran, mamini@modares.ac.ir \newline \indent Department of
Mathematics, Teacher Training University, Tehran, Iran \linebreak
medghalchi@saba.tmu.ac.ir} 
\thanks{This research was supported by Grant 510-2090
of Shahid Behshti University}

\begin{abstract}
The Fourier and Fourier-Stieltjes algebras $A(G)$ and $B(G)$ of a
locally  compact group $G$ are introduced and studied in 60's by
Piere Eymard in his PhD thesis. If $G$ is a locally compact
abelian group, then $A(G)\simeq L^1(\hat{G})$, and $B(G)\simeq
M(\hat{G})$, via the Fourier and Fourier-Stieltjes transforms,
where $\hat{G}$ is the Pontryagin dual of $G$. Recently these
algebras are defined on a (topological or measured) groupoid and
have shown to share many common features with the group case. This
is the last in a series of papers in which we have investigated a
"restricted" form of these algebras on a unital inverse semigroup
$S$.
\end{abstract}

\maketitle

\section{Introduction.}
In [1] and [2] we introduced the concept of restricted \repns on
an inverse semigroup and studied the restricted versions of
positive definite functions, semigroup algebra, and semigroup
$C^*$-algebras.

 In this paper, our aim is to study the restricted Fourier and
Fourier-Stieltjes algebras $A(S)$ and $B(S)$ on an inverse
semigroup $S$. In particular, we prove restricted version of the
Eymard's characterization [5] of the Fourier algebra (Theorem
2.1).

The structure of algebras $B(S)$ and $A(S)$ is far from being well
understood, even in special cases. From the results of [4], [7],
it is known that for a commutative unital discrete $*$-semigroup
$S$, $B(S)=M(\hat{S})^{\hat{}}$ via Bochner theorem [7]. Even in
this case, the structure of $A(S)$ seems to be much more
complicated than the group case. This is mainly because of the
lack of an appropriate analog of the group algebra. If $S$ is a
discrete idempotent semigroup with identical involution. Then
$\hat{S}$ is a compact topological semigroup with pointwise
multiplication. We believe that in this case
$A(S)=L(\hat{S})^{\hat{}}$ where $L(\hat{S})$ is the Baker algebra
on $\hat S$ (see for instance [8]), however we are not able to
prove it at this stage.

All over this paper, $S$ denotes a unital inverse semigroup with
identity $1$. $E$ denotes the set of idempotents of $S$  consists
of elements the form $ss^*,\, s\in S$. $\Sigma=\Sigma(S)$ is the
family of all $*$-representations $\pi$ of $S$ with
$$\|\pi\|:=\sup_{x\in S}\|\pi(x)\|\leq 1.$$
The {\bf associated groupoid} of $S$ is denoted by $S_a$ [9]. If
we adjoin a zero element $0$ to this groupoid, and put $0^*=0$, we
get an inverse semigroup $S_r$ which is called the {\bf restricted
semigroup} of $S$. A {\bf restricted representation}
$\{\pi,\mathcal H_\pi\}$ of $S$ is a map $\pi:S\to\mathcal
B(\mathcal H_\pi)$ such that $\pi(x^*)=\pi(x)^*\quad (x\in S)$ and
$$\pi(x)\pi(y)=\begin{cases}
\pi(xy) & \text{if}\;\; x^*x=yy^* \\
0 & \text{otherwise}
\end{cases}\quad (x,y\in S).$$

$\Sigma_r=\Sigma_r(S)$ denotes the family of all restricted
representations $\pi$ of $S$ with $\|\pi\|=\sup_{x\in
S}\|\pi(x)\|\leq 1$. Two basic examples of restricted
representations are the restricted left and right regular
representations $\lambda_r$ and $\rho_r$ of $S$ [1].

\section{Restricted Fourier and Fourier-Stieltjes algebras}

Let $S$ be a unital inverse semigroup and $P(S)$ be the set of all
bounded positive definite functions on $S$ (see [6] for the group
case and [2] for invese semigroups). Following the notations of
[2], we use the notation $P(S)$ with indices $r$, $e$, $f$ , and
$0$ to denote the positive definite functions which are
restricted, extendible, of finite support, or vanishing at zero,
respectively. Let $B(S)$ be the linear span of $P(S)$. Then $B(S)$
is a commutative Banach algebra with respect to the pointwise
multiplication and the following norm [3],[13]
$$\|u\|=\sup\{|\sum_{x\in S} u(x)f(x)| : f\in
\ell ^1(S), \sup_{\pi\in\Sigma(S)} \|\tilde{\pi}(f)\|\leq 1\}
\quad (u\in B(S)).$$ Also $B(S)$ coincides with the set of the
coefficient functions of elements of $\Sigma(S)$ [1]. If one wants
to get a similar result for the set of coefficient functions of
elements of $\Sigma_r(S)$, one has to apply the above facts to
$S_r$. But $S_r$ is not unital in general, so one is led to
consider a smaller class of bounded positive definite functions on
$S_r$. The results of [10] suggests that these should be the class
of extendible positive definite functions on $S$. Among these,
those which vanish at $0$ correspond to elements of $P_{r,e}(S)$.

In this section we show that the linear span $B_{r,e}(S)$ of
$P_{r,e}(S)$ is a commutative Banach algebra with respect to the
pointwise multiplication and an appropriate modification of the
above norm. We call this the {\bf restricted Fourier-Stieltjes
algebra} of $S$ and show that it coincides with the set of all
coefficient functions of elements of $\Sigma_r(S)$.

As before, the indices $e$, $0$, and $f$ is used to distinguish
extendible elements, elements  vanishing at $0$, and elements of
finite support, respectively. We freely use any combination of
these indices.  Consider the linear span of $P_{r,e,f}(S)$ which
is clearly a two-sided ideal of $B_{r,e}(S)$, whose closure
$A_{r,e}(S)$ is called the {\bf restricted Fourier algebra} of
$S$. We show that it is a commutative Banach algebra under
pointwise multiplication and  norm of $B_r(S)$. We also show that
it is the predual of the von Neumann algebra.

In order to study properties of $B_{r,e}(S)$, we are led by
Proposition 5.1 to consider $B_{0,e}(S_r)$. More generally we
calculate this algebra for any inverse $0$-semigroup $T$. Let
$B_e(T)$ be the linear span of $P_e(T)$ with pointwise
multiplication and the norm
$$\|u\|=sup\{\,\big|  \sum_{x\in T} f(x)u(x)\big|: f\in\ell^1(T),
\|f\|_{\Sigma(T)}\leq 1\}\quad(u\in B_e(T)),$$
and $B_{0,e}(T)$ be the closed ideal of $B_e(T)$ consisting of
elements vanishing at $0$. First let us show that $B_e(T)$ is
complete in this norm. The next lemma is quite well known and
follows directly from the definition of the functional norm.

\begin{lemm}
If $X$ is a Banach space, $D\subseteq X$ is dense, and $f\in X^*$,
then
$$\|f\|=sup\{|f(x)|: x\in D,\|x\|\leq 1\}.$$\qed
\end{lemm}

\begin{lemm}
If $T$ is an inverse $0$-semigroup (not necessarily unital), then
we have the following isometric isomorphism of Banach spaces:

$(i)$ $B_e(T)\simeq C^*(T)^*$,

$(ii)$ $B_{0,e}(T)\simeq \big(C^*(T)/\mathbb C\delta_0\big)^*$.

In particular $B_e(T)$ and $B_{0,e}(T)$ are Banach spaces.
\end{lemm}
{\bf Proof} $(ii)$ clearly follows from $(i)$. To prove $(i)$,
first recall that $P_e(S)$ is affinely isomorphic to
$\ell^1(S)_{+}^*$ [10, 1.1] via
$$<u,f>=\sum_{x\in S} f(x)u(x)\quad(f\in\ell^1(S), u\in P_e(S)).$$
This defines an isometric isomorphism $\tau_0$ from $B_e(T)$ into
$\ell^1(T)^*$ (with the dual norm). By above lemma, one can lift
$\tau_0$ to an isometric isomorphism $\tau$ from $B_e(T)$ into
$C^*(T)^*$. We only need to check that $\tau$ is surjective. Take
any $v\in C^*(T)$, and let $w$ be the restriction of $v$ to
$\ell^1(T)$. Since $\|f\|_{\Sigma(T)}\leq\|f\|_1$, for each
$f\in\ell^1(T)$, The norm of $w$ as a linear functional on
$\ell^1(T)$ is not bigger than of the norm of $v$ as a functional
on $C^*(T)$. In particular, $w\in\ell^1(T)^*$ and so there is a
$u\in B_e(T)$ with $\tau_0(u)=w$. Then $\tau(u)=v$, as
required.\qed

According to Notation 2.1 of [2], we know that the restriction map
$\tau:B_{0,e}(S_r)\to B_{r,e}(S)$ is a surjective linear
isomorphism. Also $\tau$ is clearly an algebra homomorphism
($B_{0,e}(S_r)$ is an algebra under pointwise multiplication [10,
3.4], and the surjectivity of $\tau$ implies that the same fact
holds for $B_{r,e}(S)$ \,). Now we put the following norm on
$B_{r}(S)$,
$$\|u\|_r=sup\{\,\big|\sum_{x\in S} f(x)u(x)\big|: f\in\ell_r^1(S),
\|f\|_{\Sigma_r(S)}\}\quad (u\in B_r(S)),$$ then using the fact
that $B_{0,e}(S_r)$ is a Banach algebra (it is a closed subalgebra
of $B(S_r)$ which is a Banach algebra [3, Theorem 3.4]) we have

\begin{lemm} The restriction map $\tau:B_{0,e}(S_r)\to B_{r,e}(S)$ is
an isometric isomorphism of normed algebras. In particular, $B_{r,e}(S)$
is a commutative Banach algebra under pointwise multiplication and above norm.
\end{lemm}
{\bf Proof} The second assertion follows from the first and Lemma
2.2 applied to $T=S_r$. For the first assertion, we only need to
check that $\tau$ is an isometry. But this follows directly from
[1, Proposition 3.3] and the fact that
$\Sigma_r(S)=\Sigma_0(S_r)$.\qed

\begin{cor}
$B_{r,e}(S)$ is the set of coefficient functions of elements of
$\Sigma_r(S)$.
\end{cor}
{\bf Proof} Given $u\in P_{r,e}(S)$, let $v$ be the extension by
zero of $u$ to a function on $S_r$, then $v\in P_{0,e}(S_r)$, so
there is a cyclic \repn $\pi\in\Sigma(S_r)$, say with cyclic
vector $\xi\in\mathcal H_\pi$, such that $v=<\pi(.)\xi,\xi>$ (see
the proof of [10, 3.2]). But
$$0=v(0)=<\pi(0)\xi,\xi>=<\pi(0^*0)\xi,\xi>=\|\pi(0)\xi\|,$$ that
is $\pi(0)\xi=0$. But $\xi$ is the cyclic vector of $\pi$, which
means that for each $\eta\in\mathcal H_\pi$, there is a net of
elements of the form $\sum_{i=1}^n c_i\pi(x_i)\xi$, converging to
$\eta$ in the norm topology of $\mathcal H_\pi$, and
$$\pi(0)(\sum_{i=1}^n c_i\pi(x_i)\xi=\sum_{i=1}^n c_i\pi(0)\xi=0,$$
so $\pi(0)\eta=0$, and so $\pi(0)=0$. This means that
$\pi\in\Sigma_0(S_r)=\Sigma_r(S)$. Now a standard argument, based
on the fact that $\Sigma_r(S)=\Sigma_0(S_r)$ is closed under
direct sum, shows that each $u\in B_{r,e}(S)$ is a coefficient
function of some element of $\Sigma_r(S)$. The converse follows
from [2, Lemma 2.6].\qed

\begin{cor}
We have the isometric isomorphism of Banach spaces
$B_{r,e}(S)\simeq C_r(S)^*$.
\end{cor}
{\bf Proof} We have the following of isometric linear
isomorphisms: first $B_{r,e}(S)\simeq B_{0,e}(S_r)$ (Lemma 2.3),
then $B_{0,e}(S_r)\simeq \big(C^*(S_r)/\mathbb C\delta_0\big)^*$
(Lemma 2.2, applied to $T=S_r$), and finally $C_r^*(S)\simeq
C^*(S_r)/\mathbb C\delta_0$ [1, Proposition 4.2].\qed

Next, as in [7], we give an alternative description of the norm of
the Banach algebra $B_{r,e}(S)$. For this we need to know more
about the universal \repn of $S_r$. Applying the discussion before
Example 2.1 in [2] to $T=S_r$, we know that the universal \repn
$\omega$ of $S_r$ is the direct sum of all cyclic \repns
corresponding to elements of $P_e(S_r)$. To be more precise, this
means that given any $u\in P_e(S_r)$ we consider $u$ as a positive
linear functional on $\ell^1(S_r)$, then by [10, 21.24], there is
a cyclic \repn $\{\tilde \pi_u,\mathcal H_u,\xi_u\}$ of
$\ell^1(S_r)$, with $\pi_u\in\Sigma(S_r)$, such that
$$<u,f>=<\tilde\pi_u(f)\xi_u,\xi_u>\quad(f\in\ell^1(S_r)).$$
Therefore $\pi_u$ is a cyclic \repn of $S_r$ and
$u=<\pi_u(.)\xi_u,\xi_u>$ on $S_r$. Now $\omega$ is the direct sum
of all $\pi_u$'s, where $u$ ranges over $P_e(S_r)$. There is an
alternative construction in which one can take the direct sum of
$\pi_u$'s with $u$ ranging over $P_{0,e}(S_r)$ to get a sub\repn
$\omega_0$ of $\omega$. Clearly $\omega\in\Sigma(S_r)$ and
$\omega_0\in\Sigma_0(S_r)$. It follows from [10,3.2] that the set
of coefficient functions of $\omega$ and $\omega_0$ are $B_e(S_r)$
and $B_{0,e}(S_r)=B_{r,e}(S)$, respectively (c.f. Notation 2.1 in
[2]) . As far as the original semigroup $S$ is concerned, we
prefer to work with $\omega_0$, since it could be considered as an
element of $\Sigma_r(S)$. Now $\tilde\omega_0$ is a non degenerate
$*$-\repn of $\ell_r^1(S)$ which uniquely extends to a non
degenerate \repn of the restricted full $C^*$-algebra $C_r^*(S)$,
which we still denote by $\tilde\omega_0$. We gather some of the
elementary facts about  $\tilde\omega_0$ in the next lemma.

\begin{lemm} With the above notation, we have the following:

$(i)$ $\tilde\omega_0$ is the direct sum of all non degenerate
\repns $\pi_u$ of $C_r^*(S)$ associated with elements $u\in
C_r^*(S)_{+}^*$ via the GNS-construction, namely $\tilde\omega_0$
is the universal \repn of $C_r^*(S)$. In particular, $C_r^*(S)$ is
faithfully represented in $\mathcal H_{\omega_0}$.

$(ii)$ The von Numann algebras $C_r^*(S)^{**}$  and the double
commutant of $C_r^*(S)$ in $\mathcal B(\mathcal H_{\omega_0})$ are
isomorphic. They are generated by elements $\tilde\omega_0(f)$,
wuth $f\in\ell_r^1(S)$, as well as by elements $\omega_0(x)$, with
$x\in S$.

$(iii)$ Each \repn $\pi$ of $C_r^*(S)$ uniquely decomposes as
$\pi=\pi^{**}\circ\omega_0$.

$(iv)$ For each $\pi\in\Sigma_r(S)$ and $\xi,\eta\in\mathcal
H_\pi$, let $u=<\pi(.)\xi,\eta>$, then $u\in C_r(S)^*$ and $$<T,u>
= <\tilde\pi^{**}(T)\xi,\eta> \quad (T\in C_r^*(S)^{**}).$$
\end{lemm}
{\bf Proof} Statement $(i)$ follows by an standard argument.
Statement $(iii)$ and The first part of $(ii)$ follow from $(i)$
and the second part of $(ii)$ follows from the fact that both set
of elements described in $(ii)$ have clearly the same commutant in
$\mathcal B(\mathcal H_{\omega_0})$ as the set of elements
$\tilde\omega_0(u)$, with $u\in C_r^*(S)$ which generate
$C_r^*(S)^{''}$. The first statement of $(iv)$ follows from [2,
Lemma 2.6] and Corollary 2.2. As for the second statement, first
note that for each $f\in\ell_r^1(S)$, $\tilde\omega_0(f)$ is the
image of $f$ under the canonical embedding of $C_r^*(S)$ in
$C_r^*(S)^{**}$. Therefore, by $(iii)$,
\begin{align*}
<\tilde\omega_0(f),u>&=<u,f>=\sum_{x\in S} f(x)u(x)\\
&=<\tilde\pi(f)\xi,\eta>=<\tilde\pi^{**}\circ\tilde\omega_0(f)\xi,\eta>.
\end{align*}
Taking limit in $\|.\|_{\Sigma_r}$ we get the same relation for
any $f\in C_r^*(S)$, and then, using $(ii)$, by taking limit in
the ultraweak topology of $C_r^*(S)^{**}$, we get the desired
relation.\qed

\begin{lemm} Let $1$ be the identity of $S$, then for each $u\in P_{r,e}(S)$
we have $\|u\|_r=u(1)$.
\end{lemm}
{\bf Proof} As $\|\delta_e\|_{\Sigma_r}=1$ and
$u(1)=\lambda_r(1)u(1)\geq 0$, we have $\|u\|_r\geq |u(1)|=u(1)$.
Conversely, by the proof of Corollary 2.1, there is
$\pi\in\Sigma_r(S)$ and $\xi\in\mathcal H_\pi$ such that
$u=<\pi(.)\xi,\xi>$. Hence $u(1)=<\pi(1)\xi,\xi>=\|\xi\|^2 \geq
\|u\|_r$.\qed

\begin{lemm}
For each $\pi\in\Sigma_r(S)$ and $\xi,\eta\in\mathcal H_\pi$,
consider $u=<\pi(.)\xi,\eta>\in B_{r,e}(S)$, then $\|u\|_r\leq
\|\xi\|.\|\eta\|$ . Conversely each $u\in B_{r,e}(S)$ is of this
form and and we may always choose $\xi,\eta$ so that
$\|u\|_r=\|\xi\|.\|\eta\|$.
\end{lemm}
{\bf Proof} The first assertion follows directly from the
definition of $\|u\|_r$ (see the paragraph after Lemma 2.2). The
first part of the second assertion is the content of Corollary
2.1. As for the second part, basically the proof goes as in [5].
Consider $u$ as an element of $C_r^*(S)^*$ and let $u=v.|u|$ be
the polar decomposition  of $u$, with $v\in C_r^*(S)^{**}$ and
$|u|\in C_r^*(S)_{+}^*=P_{r,e}(S)$, and the dot product is the
module action of $C_r^*(S)^{**}$ on $C_r^*(S)^*$ [5]. Again, by
the proof of Corollary 2.1, there is a cyclic \repn $\pi\in
\Sigma_r(S)$, say with cyclic vector $\eta$, such that
$|u|=<\pi(.)\eta,\eta>$. Put $\xi=\tilde\pi^{**}(v)\eta$, then
$\|\xi\|\leq\|\eta\|$ and by Lemma 2.4 $(iv)$ applied to $|u|$,
$$u(x)=<\omega_0(x),u>=<\omega_0(x)v,|u|>=<\tilde\pi^{**}
\circ\omega_0(x)(v)\eta,\eta>=<\pi(x)\xi,\eta>,$$ and, by
Corollary 2.2 and Lemma 2.5,
$$\|u\|_r=\|\,|u|\,\|=|u|(1)=\|\eta\|^2\geq \|\xi\|.\|\eta\|.$$\qed

Note that the above lemma provides an alternative (direct) way of
proving the second statement of Lemma 2.3 (just take any two
elements $u,v$ in $B_{r,e}(S)$ and represent them as coefficient
functions of two \repns such that the equality hold for the norms
of both $u$ and $v$, then use the tensor product of those \repns
to represent $uv$ and apply the first part of the lemma to $uv$.)
Also it gives the alternative description of the norm on
$B_{r,e}(S)$ as follows:

\begin{cor} For each $u\in B_{r,e}(S)$,
$$\|u\|_r=inf\{\|\xi\|.\|\eta\|: \xi,\eta\in\mathcal H_\pi, \,\pi\in \Sigma_r(S),
\, u=<\pi(.)\xi,\eta>\}.$$\qed
\end{cor}

\begin{cor} For each $u\in B_{r,e}(S)$,
$$\|u\|_r=sup \{\big|\sum_{n=1}^\infty c_n u(x_n)\big|: c_n\in\mathbb C, x_n\in S
\,(n\geq 1), \|\sum_n c_n\delta_{x_n}\|_{\Sigma_r}\leq 1\}.$$
\end{cor}
{\bf Proof} Just apply Kaplansky's density theorem to the unit
ball of $C_r^*(S)^{**}$.\qed

\begin{cor}
The unit ball of $B_{r,e}(S)$ is closed in the topology of
pointwise convergence.
\end{cor}
{\bf Proof} If $u\in B_{r,e}(S)$ with $\|u\|_r\leq 1$, then for
each $n\geq 1$, each $c_1,\dots,c_n\in\mathbb C$, and each
$x_1,\dots,x_n\in S$,
$$\big|\sum_{k=1}^n c_k u(x_k)\big|\leq \big\|\sum_{k=1}^n c_k\delta_{x_k}\big\|_{\Sigma_r}.$$
If $u_\alpha\to u$, pointwise on $S$ with $u_\alpha\in
B_{r,e}(S)$, $\|u_\alpha\|_r\leq 1$, for each $\alpha$, then all
$u_\alpha$'s satisfy above inequality, and so does $u$. Hence, by
above corollary, $u\in B_{r,e}(S)$ and $\|u\|_r\leq 1$.\qed

\begin{lemm}
For each $f,g\in \ell^2(S)$, $f\bullet\tilde{g}\in B_{r,e}(S)$ and
if $\|\cdot\|_r$ is the norm of $B_{r,e}(S)$,
$\|f\bullet\tilde{g}\|_r\leq\|f\|_2 .\|g\|_2.$
\end{lemm}
{\bf Proof} The first assertion follows from polarization identity
[2, Lemma 5.3] and the fact that for each $h\in \ell^2(S)$,
$h\bullet\tilde{h}$ is a restricted extendible positive definite
function [2, Theorem 5.1]. Now if $u=f\bullet\tilde{g}$, then
\begin{align*}
\|u\|_r &= sup\{\, \big|\sum_{y\in S} u(y)\varphi(y)\big|:
\varphi\in\ell_r^1(S),
\|\varphi\|_{\Sigma_r}\leq 1\}\\
&= sup\{\,\big|\sum_{y\in S}
<\lambda_r(y^*)f,g>\varphi(y)\big|:\varphi\in\ell_r^1(S),
\|\varphi\|_{\Sigma_r}\leq 1\}\\
& = sup\{\,\big|\sum_{y\in S} <f,\lambda_r(y)g>\varphi(y)\big|:
\varphi\in\ell_r^1(S),
\|\varphi\|_{\Sigma_r}\leq 1\}\\
& = sup\{\,|<f,\tilde\lambda_r(\varphi)g>|:
\varphi\in\ell_r^1(S),
\|\varphi\|_{\Sigma_r}\leq 1\}\\
& = sup_{\|\varphi\|_{\Sigma_r}\leq 1}
\|\tilde\lambda_r(\varphi)\|\|f\|_2\|g\|_2\leq
\|f\|_2.\|g\|_2.\qed
\end{align*}
The next theorem extends Eymard's theorem [5, 3.4] to inverse
semigroups.

\begin{theo}
Consider the following sets:
\begin{align*}
E_1 & =<f\bullet\tilde{g}: f,g\in \ell_{f}^2(S)>, \\
E_2 & = <h\bullet\tilde{h}: h\in \ell_{f}^2(S)>, \\
E_3 & = <P_{r,e,f}(S)>,\\
E_4 & = <P(S)\cap \ell^2(S)>\\
E_5 & =<h\bullet\tilde{h}:h\in \ell^2(S)>\\
E_6 & =<f\bullet\tilde{g}:f,g\in \ell^2(S)>
\end{align*}
Then $E_1\subseteq E_2\subseteq E_3\subseteq E_4\subseteq
E_5\subseteq E_6\subseteq B_{r,e}(S)$ and the closures of all of
these sets in $B_{r,e}(S)$ are equal to $A_{r,e}(S)$.
\end{theo}
{\bf Proof} The inclusion $E_1\subseteq E_2$ follows from [2,
Lemma 2.3], and the inclusions $E_2\subseteq E_3$ and
$E_4\subseteq E_5$ follow from [2, Theorem 2.1]. The inclusions
$E_3\subseteq E_4$ and $E_5\subseteq E_6$ are trivial. Now $E_1$
is dense in $E_6$ by lemma 2.7 and the fact that $\ell_f^2(S)$ is
dense in $\ell^2(S)$. Finally $\bar{E}_3=A_{r,e}(S)$, by
definition, and $E_3\subseteq E_2\subseteq E_1$ by [2, Theorem
2.1], hence $\bar{E}_i=A_{r,e}(S)$, for each $1\leq i\leq 6$. \qed

\begin{lemm} $P_{r,e}(S)$ separates the points of $S$.
\end{lemm}
{\bf Proof} We know that $S_r$ has a faithful \repn (namely the
left regular \repn $\Lambda$), so $P_e(S_r)$ separates the points
of $S_r$ [10, 3.3]. Hence $P_{0,e}(S_r)=P_{r,e}(S)$ separates the
points of $S_r\backslash\{0\}=S$.\qed

\begin{prop}
For each $x\in S$ there is $u\in A_{r,e}(S)$ with $u(x)=1$. Also
$A_{r,e}(S)$ separates the points of $S$.
\end{prop}
{\bf Proof} Given $x\in S$, let
$u=\delta_{(x^*x)}\bullet\tilde\delta_{x^*}\in E_1\subseteq
A_{r,e}(S)$, then $u(x)=1$. Also given $y\neq x$ and $u$ as above,
if $u(y)\neq 1$, then $u$  separates $x$ and $y$. If $u(y)=1$,
then use above lemma to get some $v\in B_{r,e}(S)$ which separates
$x$ and $y$. Then $u(x)=u(y)=1$, so $(uv)(x)=v(x)\neq
v(y)=(uv)(y)$, i.e. $uv\in A_{r,e}(S)$ separates $x$ and $y$. \qed

\begin{prop}
For each finite subset $K\subseteq S$, there is $u\in
P_{r,e,f}(S)$ such that $u|_K\equiv 1$.
\end{prop}
{\bf Proof} For $F\subseteq S$, let $F_e=\{x^*x: x\in F\}$ and
note that $F\subseteq F\bullet F_e$ (since $x=x(x^*x)$, for each
$x\in F$). Now given a finite set $K\subseteq S$, put $F=K\cup
K^*\cup K_e$, then since $K_e=K_e^*$ we have $F=F^*$, and since
$K_e=(K^*)_e$ and $(K_e)_e=K_e$ we have $F_e\subseteq F$. Hence
$K\subseteq F\subseteq F\bullet F$. Now $F\bullet F$ is a finite
set and if $f=\chi_F$, then $u=f\bullet
\tilde{f}=\chi_F\bullet\tilde{\chi}_F=\chi_{F\bullet
F^*}=\chi_{F\bullet F}\in P_{r,e,f}(S)$ and $u|_K\equiv 1$. \qed

\begin{cor}  $B_{r,e,f}(S)=<P_{r,e,f}(S)>$\,
and $\overline{B_{r,f}(S)}$ $=A_{r,e}(S)$.
\end{cor}
{\bf Proof} Clearly $<P_{r,e,f}(S)>\subseteq B_{r,e,f}(S)$. Now if
$v\in B_{r,e,f}(S)$, then $v=\sum_{i=1}^4 \alpha_iv_i$, for some
$\alpha_i\in\BC$ and $v_i\in P_{r,e,f}(S)$ $(1\leq i\leq 4)$. Let
$K=\supp(v)\subseteq S$ and $u\in P_{r,e,f}(S)$ be as in the above
proposition, then $u|_K\equiv 1$ so
$v=uv=\sum_{i=1}^4\alpha_i(uv_i)$ is in the linear span of
$P_{r,e,f}(S)$. \qed

\section{Fourier and Fourier-Stieltjes algebras of associated groupoids}

We observed in section 1 that one can naturally associate a
(discrete) groupoid $S_a$ to any inverse semigroup $S$. The
Fourier and Fourier-Stieltjes algebras of (topological and
measured) groupoids are studied in [11], [12], [13], and [14]. It
is natural to ask if the results of these papers, applied to the
associated groupoid $S_a$ of $S$, could give us some information
about the associated algebras on $S$. In this section we explore
the relation between $S$ and its associated groupoid $S_a$, and
resolve some technical difficulties which could arise when one
tries to relate the corresponding function algebras.

Let us remind some general terminology and facts about groupoids.
There are two parallel approaches to the theory of groupoids,
theory of measured groupoids versus theory of locally compact
groupoids (compare [13] with [14]). Here we deal with discrete
groupoids (like $S_a$) and so basically it doesn't matter which
approach we take, but the topological approach is more suitable
here. Even if one wants to look at the topological approach, there
are two different interpretation about what we mean by a "\repn"
(compare [12] with [13]). The basic difference is that whether we
want \repns to preserve multiplications {\it everywhere} or just
{\it almost everywhere} (with respect to a Borel measure on the
unit space of our groupoid which changes with each \repn ). Again
the "everywhere approach" is more suitable for our setting. This
approach, mainly taken by [11] and [12], is the best fit for the
\repn theory of inverse semigroups (when one wants to compare
\repn theories of $S$ and $S_a$). Even then, there are some basic
differences which one needs to deal with them carefully.

We mainly follow the approach and terminology of [12]. As we only
deal with discrete groupoids we drop the topological
considerations of [12]. This would simplify our short introduction
and facilitates our comparison. A (discrete) {\bf groupoid} is a
set $G$ with a distinguished subset $G^2\subseteq G\times G$ of
pairs of multiplicable elements, a multiplication map $:G^2\to G ;
\, (x,y)\mapsto xy$, and an inverse map $:G\to G ; \, x\mapsto
x^{-1}$\,, such that for each $x,y,z\in G$

$(i)$ $(x^{-1})^{-1}=x$,

$(ii)$ If $(x,y),(y,z)$ are in $G^2$, then so are $(xy,z),
(x,yz)$, and $(xy)z=x(yz)$,

$(iii)$ $(x^{-1},x)$ is in $G^2$ and if $(x,y)$ is in $G^2$ then
$x^{-1}(xy)=y$,

$(iv)$ If $(y,x)$ is in $G^2$ then $(yx)x^{-1}=y$.

For $x\in G$, $s(x)=x^{-1}x$ and $r(x)=xx^{-1}$ are called the
source and range of $x$, respectively. $G^0=s(G)=r(G)$ is called
the unit space of $G$. For each $u,v\in G^0$ we put
$G^u=r^{-1}(u)$, $G_v=s^{-1}(v)$, and $G_v^u=G_v\cap G^u$. Note
that for each $u\in G^0$, $G_u^u$ is a (discrete) group, called
the isotropy group at $u$. Any (discrete) groupoid $G$ is endowed
with left and right Haar systems $\{\lambda_u\}$ and
$\{\lambda^u\}$, where $\lambda_u$ and $\lambda^u$ are simply
counting measures on $G_u$ and $G^u$, respectively. Consider the
algebra $c_{00}(G)$ of finitely supported functions on $G$. We
usually make this into a normed algebra using the so-called
$I$-norm
$$\|f\|_I=max\{sup_{u\in G^0}\sum_{x\in G_u} |f(x)|, \, sup_{u\in G^0}
\sum_{x\in G^u} |f(x)|\}\quad (f\in c_{00}(G)),$$
where the above supremums are denoted respectively by
$\|f\|_{I,s}$ and $\|f\|_{I,r}$. Note that in general $c_{00}(G)$
is not complete in this norm. We show the completion of
$c_{00}(G)$ in $\|.\|_I$ by $\ell^1(G)$. There are also natural
$C^*$-norms in which one can complete $c_{00}(G)$ and get a
$C^*$-algebra. Two well known groupoid $C^*$-algebras obtained in
this way are the {\bf full and reduced groupoid} $C^*$-algebras
$C^*(G)$ and $C^*_L(G)$. Here we briefly discuss their
construction and refer the reader to [] for more details.

A Hilbert bundle $\mathcal H=\{\mathcal H_u\}$ over $G^0$  is just
a field of Hilbert spaces indexed by $G^0$. A {\bf \repn} of $G$
is a pair $\{\pi,\mathcal H^\pi\}$ consisting of a map $\pi$ and a
Hilbert bundle $\mathcal H^\pi=\{\mathcal H_u^\pi\}$ over $G^0$
such that For each $x,y\in G$,

$(i)$  $\pi(x):\mathcal H^\pi_{s(x)}\to\mathcal H^\pi_{r(x)}$ is a
surjective linear isometry,

$(ii)$ $\pi(x^{-1})=\pi(x)^*$,

$(iii)$ If $(x,y)$ is in $G^2$, $\pi(xy)=\pi(x)\pi(y)$.

We usually just refer to $\pi$ as the \repn and it is always
understood that there is a Hilbert bundle involved. We denote the
set of all \repns of $G$ by $\Sigma(G)$. Note that here a \repn
corresponds to a (continuous) Hilbert bundle, where as in the
usual approach to (locally compact or measured) categories \repns
are given by measurable Hilbert bundles (see [12] for more
details).

A natural example of such a \repn is the {\bf left regular \repn}
$L$ of $G$. The Hilbert bundle of this \repn is $L^2(G)$ whose
fiber at $u\in G^0$ is $L^2(G^u,\lambda^u)$. In our case that $G$
is discrete, this is simply $\ell^2(G^u)$. Each $f\in c_{00}(G)$
could be regarded as a section of this bundle (which sends $u\in
G^0$ to the restriction of $f$ to $G^u$). Also $G$ acts on bounded
sections $\xi$ of $L^2(G)$ via
$$L_x\xi(y)=\xi(x^{-1}y)\quad (x\in G, y\in G^{r(x)}).$$
Let $E^2(G)$ be the set of sections of $L^2(G)$ vanishing at
infinity. This is a Banach space under the sup-norm and contains
$c_{00}(G)$. Furthermore, it is a canonical $c_0(G^0)$-module via
$$b\xi(x)=\xi(x)b(r(x))\quad (x\in G, \xi\in E^2(G), b\in c_0(G^0)).$$
Now $E^2(G)$ with the $c_0(G)$-valued inner product
$$<\xi,\eta>(u)=<L(.)\xi^u\circ s(.),\eta^u\circ r(.)>$$
is a Hilbert $C^*$-module. The action of $c_{00}(G)$ on itself by
left convolution extends to a $*$-anti\repn of $c_{00}(G)$ in
$E^2(G)$, which is called the left regular \repn of $c_{00}(G)$
[12, Proposition 10]. The map $f\mapsto L_f$ is a norm decreasing
homomorphism from $(c_{00}(G), \|.\|_{I,r})$ into $\mathcal
B(E^2(G))$. Also the former has a left bounded approximate
identity $\{e_\alpha\}$ consisting of positive functions such that
$\{L_{e_\alpha}\}$ tends to the identity operator in the strong
operator topology of the later [12, Proposition 11]. The closure
of the image of $c_{00}(G)$ under $L$ is a $C^*$-subalgebra
$C^*_L(G)$ of $\mathcal B(E^2(G))$ which is called the {\bf
reduced} $C^*$-{\bf algebra} of $G$. We should warn the reader
that $\mathcal B(E^2(G))$ is merely a $C^*$-algebra and, in
contrast with the Hilbert space case, it is not a von Neumann
algebra in general. The above construction simply means that we
have used the \repn $L$ to introduce an auxiliary $C^*$-norm on
$c_{00}(G)$ and took the completion of $c_{00}(G)$ with respect to
this norm. A similar construction using all non degenerate
$*$-\repns of $c_{00}(G)$ in Hilbert $C^*$-modules yields a
$C^*$-completion $C^*(G)$ of $c_{00}(G)$, called the {\bf full}
$C^*$-{\bf algebra} of $G$.

Next one can define positive definiteness in this context. Let
$\pi\in\Sigma(G)$, for bounded sections $\xi,\eta$ of $\mathcal
H^\pi$, the function $x\mapsto <\pi(x)\xi(s(x)),\eta(r(x))>$
(where the inner product is taken in the Hilbert space $\mathcal
H^\pi_{r(x)}$) is called a coefficient function of $\pi$. A
function $\varphi\in\ell^\infty (G)$ is called {\bf positive
definite} if for all $u\in G^0$ and all $f\in c_{00}(G)$
$$\sum_{x,y\in G^u} \varphi(y^{-1}x)f(y)\bar f(x)\geq 0,$$
or equivalently, for each $n\geq 1$, $u\in G^0$, $x_1,\dots x_n\in
G^u$, and $\alpha_1,\dots,\alpha_n\in\mathbb C$
$$\sum_{i,j=1}^n \bar\alpha_i\alpha_j \varphi(x_i^{-1}x_j)\geq 0.$$
We denote the set of all positive definite functions on $G$ by
$P(G)$. The linear span $B(G)$ of $P(G)$ is called the {\bf
Fourier-Stieltjes algebra} of $G$. It is equal to the set of all
coefficient functions of elements of $\Sigma(G)$ [12, Theorem 1].
It is a unital commutative Banach algebra [12, Theorem 2] under
pointwise operations and the norm
$\|\varphi\|=inf\|\xi\|\|\eta\|$, where the infimum is taken over
all \repns $\varphi=<\pi(.)\xi\circ s(.),\eta\circ r(.)>$. On the
other hand each $\varphi\in B(G)$ could be considered as a
completely bounded linear operator on $C^*(G)$ via
$$<\varphi, f>=\varphi .f\quad(\varphi\in B(G), f\in c_{00}(G)),$$
such that $\|\varphi\|_\infty\leq\|\varphi\|_{cb}\leq\|\varphi\|$
[12, Theorem 3]. The last two norms are equivalent on $B(G)$ (they
are equal in the group case, but it is not known if this is the
case for groupoids). Following [12] we denote $B(G)$ endowed with
$cb$-norm with $\mathcal B(G)$. This is known to be a Banach
algebra (This is basically [13, Theorem 6.1] adapted to this
framework [12, Theorem 3]).

There are four candidates for the {\bf Fourier algebra} $A(G)$.
The first is the closure of the linear span of the coefficients of
$E^2(G)$ in $B(G)$ [14], the second is the closure of $\mathcal
B(G)\cap c_{00}(G)$ in $\mathcal B(G)$ [11], the third is the
closure of the of the subalgebra generated by the coefficients of
$E^2(G)$ in $B(G)$, and the last one is the completion of the
normed space of the quotient of $E^2(G)\hat{\bigotimes}E^2(G)$ by
the kernel of $\theta$  from $E^2(G)\hat{\bigotimes}E^2(G)$ into
$c_0(G)$ induced by the bilinear map $\theta: c_{00}(G)\times
c_{00}(G)\to c_0(G)$ defined by
$$\theta(f,g)=g*\check f\quad(f,g\in c_{00}(G)).$$

These four give rise to the same algebra in the group case. We
refer the interested reader to [12] for a comparison of these
approaches. Here we adapt the third definition. Then $A(G)$ is a
Banach subalgebra of $B(G)$ and $A(G)\subseteq c_0(G)$. Moreover
if
$$VN(G)=\{T\in\mathcal B(E^2(G)): TR_f=R_fT\,(f\in c_{00}(G))\},$$
where $R$ is the right regular \repn of $c_{00}(G)$ in $E^2(G)$,
then $VN(G)$ is the strong closure of $C^*_L(G)$ in $\mathcal
B(E^2(G))$. Note that here $VN(G)$ is not necessarily a von Numann
algebra. Also the isometric isomorphism between the linear dual of
$A(G)$ and $VN(G)$ may fail to exist, unless we replace $A(G)$
with an appropriate space [12, Theorem 4].

Now we are ready to compare the function algebras on inverse
semigroup $S$ and its associated groupoid $S_a$. We would apply
the above results to $G=S_a$. First let us look at the \repn
theory of these objects. As a set, $S_r$ compared to $S_a$ has an
extra zero element. Moreover, the product of two non zero element
of $S_r$ is $0$, exactly when it is undefined in $S_a$. Hence it
is natural to expect that $\Sigma(S_a)$ is related to
$\Sigma_0(S_r)=\Sigma_r(S)$. The major difficulty to make sense of
this relation is the fact that \repns of $S_a$ are defined through
Hilbert bundles, where as restricted \repns of $S$ are defined in
Hilbert spaces. But a careful interpretation shows that these are
two sides of one coin.

\begin{lemm} $\Sigma_r(S)=\Sigma(S_a)$.
\end{lemm}
{\bf Proof} First let us show that each $\pi\in\Sigma_r(S)$ could
be regarded as an element of $\Sigma(S_a)$. Indeed, for each $x\in
S$, $\pi(x):\mathcal H_\pi\to\mathcal H_\pi$ is a partial
isometry, so if we put $\mathcal H_u=\pi(u)\mathcal H_\pi\, (u\in
E_S)$, then we could regard $\pi(x)$ as an isomorphism from
$\mathcal H_{x^*x}\to\mathcal H_{xx^*}$. Using the fact that the
unit space of $S_a$ is $S_a^0=E_S$, it is easy now to check  that
$\pi\in\Sigma(S_a)$. Conversely suppose that $\pi\in \Sigma(S_a)$,
then for each $x\in S_a$, $\pi(x):\mathcal H_{s(x)}\to\mathcal
H_{r(x)}$ is an isomorphism of Hilbert spaces. Let $\mathcal
H_\pi$ be the direct sum of all Hilbert spaces $\mathcal H_u$,
$u\in E_S$, and define $\pi(x)(\xi_u)=(\eta_v)$, where
$$\eta_v=\begin{cases}
\pi(x)\xi_{x^*x} & \text{if}\;\; v=xx^* \\
0 & \text{otherwise}
\end{cases}\quad (x\in S, v\in E_S),$$
then we claim that
$$\pi(x)\pi(y)=\begin{cases}
\pi(xy) & \text{if}\;\; x^*x=yy^* \\
0 & \text{otherwise}
\end{cases}\quad (x,y\in S).$$
First let's assume that $x^*x=yy^*$, then
$\pi(xy)(\xi_u)=(\theta_v)$, where $\theta_v=0$, except for
$v=xyy^*x^*=xx^*$, for which
$\theta_v=\pi(xy)\xi_{y^*x^*xy}=\pi(xy)\xi_{y^*y}$. On the other
hand, $\pi(y)(\xi_u)=(\eta-v)$, where $\eta_v=0$, except for
$v=yy^*$, for which $\eta_v=\pi(y)\xi_{y^*y}$, and
$\pi(x)(\eta_v)=(\zeta_w)$, with $\zeta_w=0$, except for $w=xx^*$,
for which
$\zeta_w=\pi(x)\eta_{x^*x}=\pi(x)\eta_{yy^*}=\pi(x)\pi(y)\xi_{y^*y}$.
Hence $\pi(xy)(\xi_u)=\pi(x)\pi(y)(\xi_u)$, for each
$(\xi_u)\in\mathcal H_\pi$. Next assume that $x^*x\neq yy^*$, then
the second part of the above calculation clearly shows that
$\pi(x)\pi(y)(\xi_u)=0$. This shows that $\pi$ could be considered
as an element of $\Sigma_r(S)$. Finally it is clear that these two
embeddings are inverse of each other.\qed

Next, $S_r=S_a\cup\{0\}$ as sets, and for each bounded map
$\varphi:S_r\to\mathbb C$ with $\varphi(0)=0$, it immediately
follows from the definition that $\varphi\in P(S_a)$ if and only
if $\varphi\in P_0(S_r)$. Hence by above lemma we have
\begin{prop} The Banach spaces $B_r(S)=B_0(S_r)$ and $B(S_a)$ are isometrically
isomorphic. \qed
\end{prop}

This combined with [12, Theorem 2] (applied to $G=S_a$) shows that
$B_r(S)$ is indeed a Banach algebra under pointwise multiplication
and the above linear isomorphism is also an isomorphism of Banach
algebras. By [12, Theorem 1] now we conclude that

\begin{cor} $B_r(S)$ is the set of coefficient functions of $\Sigma_r(S)$.\qed
\end{cor}

{\bf acknowledgement.} The first author would like to thank
hospitality of Professor Mahmood Khoshkam during his stay in
University of Saskatchewan, where the main part of the revision
was done.


\begin{thebibliography}{99}

\bibitem{} M. Amini, A. Medghalchi, restricted algebras on inverse semigroups I,
representation theory,
preprint,
Shahid Beheshti University, 2000.


\bibitem{} M. Amini, A. Medghalchi, restricted algebras on inverse semigroups II,
positive definite functions, preprint,
Shahid Beheshti University, 2000.


\bibitem{} M. Amini, A. Medghalchi, Fourier algebras on topological foundation $*$-semigroups,
preprint, Shahid Beheshti University, 2000.


\bibitem{a} C.F. Dunkl, D.E. Rumirez, $L^\infty$-representations of commutative
semitopological semigroups, Semigroup Forum 7 (1974) 180-199.

\bibitem{a} P. Eymard,  L'algebra de Fourier d'un groupe localement
compact, Bull. Soc. Math. France, 92 (1964) 181-236.

\bibitem{a} R. Godement,  Les fonctions de type positive et la theorie des
groupes, Trans. Amer. Math. Soc. 63 (1948) 1-84.

\bibitem{a} M. Lashkarizadeh Bami,  Bochner's theorem and the Hausdorff
moment theorem on foundation topological semigroups, Can. J. Math.
37 (1985) 785-809.

\bibitem{a} M. Lashkarizadeh Bami,  Representations of foundation
semigroups and their algebras, Can. J. Math. 37 (1985) 29-47.

\bibitem{} Mark V. Lawson, Inverse semigroups, the theory of partial symmetries,
World Scientific, Singapore, 1998.

\bibitem{} R.J. Lindahl, P.H. Maserick, Positive-definite functions on involution semigroups,
Duke Math. J. 38 (1971) 771-782.

\bibitem{} K. Oty, Fourier-Stieltjes algebras of r-discrete groupoids, J. Operator Theory,
41 (1999) 175-197.

\bibitem{} A.T. Paterson, The Fourier algebra for locally compact groupoids, preprint, 2002.

\bibitem{}
A. Ramsay, M.E. Walter, Fourier-Stieltjes algebras of locally
compact groupoids, J. Functional Analysis, 148 (1997) 314-365.

\bibitem{}
J. N. Renault, The Fourier algebra of a measured groupoid and its
multipliers, J. Functional Analysis, 145 (1997) 455-490.


\end{thebibliography}
\end{document}